\documentclass{amsart}
\usepackage[francais,english]{babel}
\usepackage{hyperref}
\usepackage{amsfonts,amsmath,amsthm,amstext}

\DeclareMathAlphabet{\mathbbold}{U}{bbold}{m}{n}
\newcommand{\zero}{\mathbbold{0}}
\newcommand{\unit}{\mathbbold{1}}
\newcommand{\val}{\operatorname{val}}
\newcommand{\mrm}[1]{\text{\rm #1}}
\newcommand{\new}[1]{{\em #1}\index{#1}}
\newcommand{\eps}{\epsilon}
\newcommand{\perm}{\operatorname{perm}}

\newcommand{\R}{\mathbb{R}}
\newcommand{\C}{\mathbb{C}}
\newcommand{\rmin}{\R_{\min}}
\newcommand{\rmax}{\R_{\max}}
\newcommand{\sA}{\mathcal{A}}
\newcommand{\sL}{\mathcal{L}}
\newcommand{\opt}{\mrm{Opt}}
\newcommand{\sat}{{\mrm{Sat}}}
\newcommand{\id}{\mrm{id}}
\newtheorem{theorem}{Theorem}
\newtheorem{theoreme}{Th\'eor\`eme}
\newtheorem{corollary}{Corollary}
\theoremstyle{remark}
\newtheorem{example}{Example}
\title[Perturbations of eigenvalues of matrix pencils]{Perturbation of eigenvalues of matrix pencils and optimal assignment problem}
\author{Marianne Akian}
\address{Marianne Akian. INRIA, Domaine de Voluceau, B.P.~105, 78153 Le Chesnay Cedex, France.}
\email{Marianne.Akian@inria.fr}
\author{Ravindra Bapat}
\address{Ravindra Bapat. Indian Statistical Institute, New Delhi, 110016, India.}
\email{rbb@isid.ac.in}
\author{St\'ephane Gaubert}
\address{St\'ephane Gaubert.
INRIA, Domaine de Voluceau,
B.P.~105, 78153 Le Chesnay Cedex, France.}
\email{Stephane.Gaubert@inria.fr}
\keywords{Perturbation theory, max-plus algebra, tropical semiring,
spectral theory, Newton-Puiseux theorem, amoeba, graphs, 
optimal assignement problem, Hungarian algorithm, 
WKB asymptotics, large deviations. \selectlanguage{french} {\em Mots cl\'es}.\/ Th\'eorie des perturbations, alg\`ebre max-plus, semi-anneau tropical, th\'eorie spectrale, th\'eor\`eme de Newton-Puiseux, amibe, graphes, probl\`eme d'affectation optimale,
algorithme Hongrois, asymptotiques WKB, grandes d\'eviations}
\subjclass[2000]{Primary 47A55. Secondary 47A75, 15A22, 05C50, 12K10}

\date{February 23, 2004}
\begin{document}
\selectlanguage{english} 
\begin{abstract}
We consider a matrix pencil whose coefficients 
depend on a positive parameter $\epsilon$, and have 
asymptotic equivalents
of the form $a\epsilon^A$ when $\epsilon$ goes to zero,
where the leading coefficient $a$
is complex, and the leading exponent $A$ is real.
We show that the asymptotic equivalent of every eigenvalue
of the pencil can be determined generically from
the asymptotic equivalents
of the coefficients of the pencil. The generic leading exponents
of the eigenvalues are the ``eigenvalues'' of a min-plus matrix pencil.
The leading coefficients of the eigenvalues are the eigenvalues of auxiliary
matrix pencils,
constructed from certain optimal assignment problems.\\[1em]
\selectlanguage{french}
\begin{center}{\bf Perturbation de valeurs propres de faisceaux matriciels}\\{\bf et probl\`eme d'affectation optimale}\end{center}\ \\[1em]
{\sc R\'esum\'e.}\/ 
Nous consid\'erons un faisceau matriciel dont les coefficients, 
d\'ependant d'un param\`etre $\epsilon$,
ont des \'equivalents asymptotiques de la forme $a\epsilon^A$, 
lorsque $\epsilon$ tend vers z\'ero par valeurs positives,
le coefficient dominant $a$ \'etant complexe et l'exposant
dominant $A$ \'etant r\'eel.
Nous montrons qu'un \'equivalent asymptotique pour chacune
des valeurs propres du faisceau peut \^etre d\'etermin\'e g\'en\'eriquement
\`a partir des \'equivalents des coefficients du faisceau.
Les exposants dominants des valeurs propres sont les 
valeurs propres d'un faisceau matriciel min-plus, et les coefficients dominants
sont les valeurs propres de faisceaux auxiliaires, construits
au moyen de probl\`emes d'affectation optimale.
\end{abstract}
\maketitle
\sloppy
\begin{center} \bf Abridged French version
\end{center}
\selectlanguage{french}
Nous consid\'erons un faisceau matriciel
$\sA_\epsilon= \sA_{\epsilon,0}+ X\sA_{\epsilon,1} +\cdots+X^d\sA_{\epsilon,d}$, o\`u $X$ est une ind\'etermin\'ee,
et o\`u pour tout $0\leq k\leq d$, $\sA_{\epsilon,k}$ est
une matrice $n\times n$ dont les coefficients, $(\sA_{\epsilon,k})_{ij}$,
sont des fonctions continues \`a valeurs complexes
d'un param\`etre positif $\epsilon$.
On s'int\'eresse au comportement asymptotique,
lorsque $\epsilon$ tend vers $0$,
des \new{valeurs propres} 
$\sL_{\epsilon}$ de $\sA_\epsilon$, qui sont
par d\'efinition les racines du polyn\^ome
$\det(\sA_\epsilon)$. 
Nous supposons que pour tout $0\leq k\leq d$, on 
a des matrices
$a_k=((a_k)_{ij})\in \C^{n\times n}$ et $A_k=((A_k)_{ij})\in (\R\cup\{+\infty\})^{n\times n}$ telles que
\(
(\sA_{\eps,k})_{ij} = (a_k)_{ij}\epsilon^{(A_k)_{ij}}+o(\epsilon^{(A_k)_{ij}})\), pour $ 1\leq i,j\leq n$.
Lorsque $(A_k)_{ij}=+\infty$, cela signifie, par convention,
que $(\sA_{\eps,k})_{ij}$ est nulle au voisinage de $0$.
Nous cherchons un \'equivalent asymptotique
de la forme $\sL_\epsilon \sim \lambda\epsilon^{\Lambda}$,
avec $\lambda\in\C\setminus\{0\}$ et $\Lambda\in \R$,
pour chaque valeur propre $\sL_\epsilon$ de $\sA_\epsilon$.

Lorsque $\sA_\epsilon=\sA_{\epsilon,0}- X\id$,
o\`u $\id$ d\'esigne la matrice identit\'e,
il s'agit du probl\`eme classique de perturbation de valeur propres
de matrices. En particulier, lorsque la matrice $\sA_{\epsilon,0}$ est
de la forme $\sA_{\epsilon,0}=\sA_{0,0}+\epsilon b$, 
pour une matrice complexe $b=(b_{ij})$,
une th\'eorie initi\'ee par Vi\v sik et Ljusternik~\cite{vishik},
et compl\'et\'ee par Lidski\u\i~\cite{lidskii}
identifie les exposants dominants $\Lambda$,
pour des valeurs g\'en\'eriques des param\`etres $b_{ij}$, 
et donne des formules explicites pour les coefficients
dominants $\lambda$, faisant intervenir les valeurs propres
de certains compl\'ements de Schur construits \`a partir
de $b$.
Voir~\cite{mbo97} pour une pr\'esentation r\'ecente.
Cependant, pour certaines valeurs de $b$,
ces compl\'ements de Schur peuvent ne pas exister.
De tels cas singuliers ont \'et\'e \'etudi\'es
en particulier dans~\cite{mbo97,maedelman}. 
Dans~\cite{abg04a}, nous avons donn\'e
une premi\`ere extension du th\'eor\`eme de Lidski\u\i~\cite{lidskii},
donnant sous certaines conditions structurelles les \'equivalents
des valeurs propres de la matrice $\sA_{\epsilon,0}$. 
Nous \'etendons ici ce r\'esultat au cas des faisceaux,
ce qui permet de r\'esoudre les cas qui demeuraient
singuliers dans~\cite{abg04a}.

Afin d'\'enoncer le r\'esultat principal,
rappelons quelques notions d'alg\`ebre min-plus.
Le \new{semi-anneau min-plus}, $\rmin$, est l'ensemble
$\R\cup\{+\infty\}$ muni de l'addition
$(a,b)\mapsto a\oplus b:=\min(a,b)$ et de la multiplication
$(a,b)\mapsto a\otimes b:=a+b$. On \'ecrira parfois
$ab$ au lieu de $a\otimes b$.
On notera $\zero:=+\infty$ et $\unit:=0$ le z\'ero et l'unit\'e de $\rmin$,
respectivement.

Nous associons au faisceau $\sA_\epsilon$
le faisceau matriciel \(A=A_0\oplus XA_1 \oplus \cdots  \oplus X^dA_d
\) \`a coefficients dans $\rmin$.
Les coefficients de $A$, $A_{ij}$, sont donc des polyn\^omes
formels \`a coefficients dans $\rmin$, en l'ind\'etermin\'ee $X$.
Nous appelons \new{polyn\^ome caract\'eristique min-plus}
le permanent \(
P_A = \perm A\).
Rappelons que pour une matrice $B=(B_{ij})$
de taille $n\times n$, \`a coefficients
dans un semi-anneau quelconque, le \new{permanent} de $B$
est d\'efini comme la somme sur toutes les permutations $\sigma$
du poids $|\sigma|_B=B_{1\sigma(1)}\cdots B_{n\sigma(n)}$.
Lorsque les coefficients de $B$ appartiennent \`a $\rmin$, 
$|\sigma|_B=B_{1\sigma(1)}+\cdots +B_{n\sigma(n)}$,
et $\perm B$ est la valeur d'une affectation
optimale dans le graphe valu\'e associ\'e \`a $B$.

Si $P$ est un polyn\^ome formel \`a coefficients dans $\rmin$,
on note $\hat{P}$ la fonction polyn\^ome associ\'ee
\`a $P$. L'application $P\mapsto \hat{P}$ est un cas
particulier de transform\'ee de Legendre-Fenchel~\cite[\S~3.3.1]{bcoq}. Cuninghame-Green et Meijer~\cite{cuning80}
ont montr\'e que la fonction polyn\^ome
$\hat P(x)$ peut se factoriser de mani\`ere
unique sous la forme 
$\hat P(x)=a(x\oplus c_1)\cdots (x\oplus c_n)$,
avec $a,c_1,\ldots,c_n\in \rmin$.
Les nombres $c_1,\ldots,c_n$ sont appel\'es \new{racines} de $P$. 
Nous noterons $\gamma_1,\ldots,\gamma_N$ les racines
du polyn\^ome caract\'eristique $P_A$, aussi appel\'ees
valeurs propres du faisceau $A$.
Les racines $\gamma_i$, et donc la fonction polyn\^ome $\hat{P}_A$,
peuvent \^etre calcul\'es en temps $O(n^4d)$
en s'inspirant de la m\'ethode de Burkard et Butkovi\v{c}~\cite{bb02}.

Pour toute matrice $B\in \rmin^{n\times n}$ telle que $\perm B\neq\zero$,
nous d\'efinissons le graphe $\opt(B)$ form\'e
des arcs participant \`a une affectation optimale:
les n\oe uds de $\opt(B)$ sont $1,\ldots,n$,
et il y a un arc de $i$ \`a $j$ s'il y a une permutation $\sigma$ 
telle que $j=\sigma(i)$ et $|\sigma|_B=\perm B$.

Nous dirons que deux vecteurs
$U,V$ de dimension $n$ \`a coefficients dans $\rmin\setminus\{\zero\}$
forment une \new{paire hongroise} relativement \`a $B$
si, quels que soient $i,j$, on a $B_{ij}\geq U_iV_j$, 
et si $U_1\cdots U_n V_1\cdots V_n=\perm B$,
les produits s'entendant dans le semi-anneau min-plus.
Ainsi, $(U,V)$ n'est autre qu'une solution optimale
du probl\`eme lin\'eaire dual du probl\`eme d'affectation
(see~\cite[2.4]{br}). En particulier,
une paire hongroise existe d\`es que $\perm B \neq \zero$, 
et elle peut \^etre trouv\'ee en temps $O(n^3)$ gr\^ace
\`a l'algorithme hongrois.
Pour toute paire hongroise $(U,V)$, on d\'efinit
le {\em graphe de saturation}, $\sat(B,U,V)$, 
qui a pour n\oe uds $1,\ldots,n$, et pour arcs
les couples $(i,j)$ tels que $B_{ij}=U_iV_{j}$.

Pour chaque racine finie $\gamma$ de $P_A$, nous
d\'efinissons les graphes $\opt_{0}(\gamma),\ldots,\opt_{d}(\gamma)$:
$\opt_{k}(\gamma)$ a pour n\oe uds $1,\ldots,n$, et a
un arc de $i$ \`a $j$ si $(i,j)\in \opt(\hat A(\gamma))$
et $(A_k)_{ij}\gamma^k=\hat A_{ij}(\gamma)$. 
Pour toute paire hongroise $(U,V)$ relativement
\`a la matrice $\hat A(\gamma)$, nous d\'efinissons
aussi les graphes $\sat_{0}(\gamma,U,V),\ldots,\sat_{d}(\gamma,U,V)$,
qui sont obtenus en rempla\c{c}ant
$\opt(\hat A(\gamma))$ par $\sat(\hat A(\gamma),U,V)$
dans la d\'efinition de $\opt_{0}(\gamma),\ldots,\opt_{d}(\gamma)$.
Enfin, si $G$ est un graphe ayant pour n\oe uds $1,\ldots,n$, 
et si $b\in \C^{n\times n}$, nous d\'efinissons la matrice $b^{G}$,
telle que $(b^{G})_{ij}=b_{ij}$ si $(i,j)\in G$,
et $(b^{G})_{ij}=0$ sinon.
Dans l'\'enonc\'e du th\'eor\`eme qui suit,
les valeurs propres sont compt\'ees avec leurs multiplicit\'es.
\begin{theoreme}
\label{th-1fr}
Soit $\gamma$ une racine finie du polyn\^ome caract\'eristique
min-plus $P_A$. Pour chaque $0\leq k\leq d$, 
notons $G^k$ le graphe \'egal a $\opt_{k}(\gamma)$
ou bien \`a $\sat_{k}(\gamma,U,V)$, pour un choix quelconque
de la paire hongroise $U,V$ relativement \`a $\hat A(\gamma)$.
Consid\'erons le faisceau auxiliaire
$a^{(\gamma)}:= a_0^{G_{0}}  + Xa_1^{G_{1}} +\cdots + 
X^da_d^{G_{d}}$.
Alors, si le faisceau $a^{(\gamma)}$
a $m_\gamma$ valeurs propres non-nulles,
$\lambda_1,\ldots,\lambda_{m_\gamma}$,
le faisceau $\sA_\epsilon$ admet ${m_\gamma}$ valeurs propres
$\sL_{\epsilon,1},\ldots,\sL_{\epsilon,{m_\gamma}}$
ayant des \'equivalents respectifs de la forme
\(
\sL_{\epsilon,i}\sim \lambda_i \epsilon^{\gamma}
\). En outre, si $0$ est une valeur propre de multiplicit\'e $m_\gamma'$
du faisceau $a^{(\gamma)}$, le faisceau $\sA_\epsilon$
admet $m_\gamma'$ valeurs propres suppl\'ementaires $\sL_\epsilon$ 
telles que $\epsilon^{-\gamma}\sL_\eps$ converge vers $0$,
et toutes les autres valeurs propres $\sL_\epsilon$
de $\sA_\epsilon$ sont telles que le module de $\epsilon^{-\gamma}\sL_\eps$
tend vers l'infini.
Enfin, pour des valeurs g\'en\'eriques des param\`etres
$(a_k)_{ij}$, ${m_\gamma}$ co\"\i ncide avec la multiplicit\'e
de la racine $\gamma$, et $m_\gamma'$ co\"\i ncide avec la somme
des multiplicit\'es des racines de $P_A$ strictement sup\'erieures
\`a $\gamma$.
\end{theoreme}
Dans la version en anglais de la pr\'esente note, nous
illustrons le Th\'eor\`eme~\ref{th-1} en raffinant
des r\'esultats de Najman~\cite{najman}.

\par\medskip\centerline{\rule{2cm}{0.2mm}}\medskip
\setcounter{section}{0}

\selectlanguage{english}
We consider a matrix pencil of the form
\[
\sA_\epsilon = \sA_{\epsilon,0}+ X\sA_{\epsilon,1} +\cdots+X^d\sA_{\epsilon,d}  \enspace ,
\]
where, for every $0\leq k\leq d$, $\sA_{\epsilon,k}$ is
a $n\times n$ matrix whose coefficients, $(\sA_{\epsilon,k})_{ij}$,
are complex valued continuous functions
of a nonnegative parameter $\epsilon$,
and $X$ is an indeterminate. 
We are interested in the asymptotic behavior,
when $\epsilon$ tends to $0$,
of the \new{eigenvalues} 
$\sL_{\epsilon}$ of $\sA_\epsilon$, which are the
roots of the polynomial $\det(\sA_\epsilon)$.

We shall assume that for every $0\leq k\leq d$, matrices
$a_k=((a_k)_{ij})\in \C^{n\times n}$ and $A_k=((A_k)_{ij})\in (\R\cup\{+\infty\})^{n\times n}$ are given, so that
\[
(\sA_{\eps,k})_{ij} = (a_k)_{ij}\epsilon^{(A_k)_{ij}} +o(\epsilon^{(A_k)_{ij}})\enspace,
\qquad \mrm{ for all } 
 1\leq i,j\leq n \enspace ,
\]
when $\epsilon$ tends to $0$.
When $(A_k)_{ij}=+\infty$, this
means by convention that $(\sA_{\eps,k})_{ij}$
is zero in a neighborhood of zero.
We look for an asymptotic equivalent
of the form $\sL_\epsilon \sim \lambda\epsilon^{\Lambda}$,
with $\lambda\in\C\setminus\{0\}$ and $\Lambda\in \R$,
for every eigenvalue $\sL_\epsilon$ of $\sA_\epsilon$.

When $\sA_\epsilon=\sA_{\epsilon,0}- X\id$,
where $\id$ denotes the identity matrix,
we recover the classical problem of perturbation of eigenvalues
of matrices. In particular, 
when $\sA_{\epsilon,0}$ is affine in $\epsilon$, i.e.\ 
$\sA_{\epsilon,0}=\sA_{0,0}+\epsilon b$, 
for some complex matrix $b=(b_{ij})$,
a theory initiated by Vi\v sik and Ljusternik~\cite{vishik}
and completed by Lidski\u\i~\cite{lidskii} identifies the leading
exponents $\Lambda$
for generic values of the parameters $b_{ij}$, 
and gives explicit formul\ae\ for the leading
coefficients $\lambda$, involving the eigenvalues
of certain Schur complements constructed from $b$.
See~\cite{mbo97} for a recent presentation.
For some values of $b$, however,
these Schur complements may not be defined.
Such singular cases have received much attention, see in particular~\cite{mbo97,maedelman}.
In~\cite{abg04a}, we generalized
the theorem of~\cite{lidskii} to the case where $\sA_{\epsilon,0}$ 
is not affine in $\epsilon$. We showed that
the leading exponents $\Lambda$ of the eigenvalues
can be computed, generically, as the ``roots'' 
of a certain min-plus polynomial constructed from $\sA_\epsilon$.
We also showed that under structural conditions
(when certain bipartite graphs constructed
from the leading exponents of the entries of $\sA_{\epsilon,0}$
have a perfect matching), the leading exponents of all the eigenvalues
can be obtained by Schur complement formul\ae\
which extend the ones of~\cite{lidskii}.
However, when these structural conditions do not hold,
the question remained open.
We show in Theorem~\ref{th-1} below that
the leading coefficients of the eigenvalues
of the pencil $\sA_\epsilon$ are obtained,
generically, as the eigenvalues of certain auxiliary
pencils (independent of $\epsilon$), constructed by solving
certain optimal assignment problems. This gives
a complete generalization of the theorem of~\cite{lidskii}.
We shall illustrate Theorem~\ref{th-1} in Corollary~\ref{cor-1}, by refining
results of Najman \cite{najman} concerning a singular perturbation
of an affine pencil.

In order to state the main result, we need
some min-plus algebraic constructions.
Recall that the \new{min-plus semiring}, $\rmin$, is the set
$\R\cup\{+\infty\}$ equipped with the addition
$(a,b)\mapsto \min(a,b)$ and the multiplication
$(a,b)\mapsto a+b$. We denote by ``$\oplus$'' the min-plus
addition, and by ``$\otimes$'' or concatenation the min-plus
multiplication. We denote by $\zero=+\infty$
and $\unit=0$ the zero and unit elements of $\rmin$, respectively.
The {\em max-plus semiring}, $\rmax$, is obtained by replacing
``$\min$'' by ``$\max$'' and $+\infty$ by $-\infty$ in this definition.

We associate to the matrix pencil $\sA_\epsilon$
the matrix pencil with coefficients
in $\rmin$,
\(
A=A_0\oplus XA_1 \oplus \cdots  \oplus X^dA_d
\).
Here, the entries of $A$, $A_{ij}$, are formal polynomials
in the indeterminate $X$
with coefficients in $\rmin$.  We consider the formal min-plus \new{characteristic polynomial}
\(
P_A = \perm A\).
Recall that for a $n\times n$ matrix $B=(B_{ij})$ 
with entries in any semiring, the \new{permanent} of $B$
is defined as the sum over all permutations $\sigma$
of the weight $|\sigma|_B=B_{1\sigma(1)}\cdots B_{n\sigma(n)}$.
Thus, if $B$ is a matrix with entries in $\rmin$, 
the weight of $\sigma$, $|\sigma|_B$,
is the usual sum $B_{1\sigma(1)}+\cdots +B_{n\sigma(n)}$,
and $\perm B$ is the value of an optimal
assignment in the weighted graph associated to $B$.
(See~\cite{br} for more background.)

If $P$ is a formal polynomial with coefficients in $\rmin$,
we denote by $\hat{P}$ the polynomial function associated
to $P$. The map $P\mapsto \hat{P}$ is a specialization
of the Legendre-Fenchel transform~\cite[\S~3.3.1]{bcoq}.
Cuninghame-Green and Meijer~\cite{cuning80} have shown that
the min-plus polynomial
function $\hat P(x)$, can be factored uniquely as 
$\hat P(x)=a(x\oplus c_1)\cdots (x\oplus c_n)$.
with $a,c_1,\ldots,c_n\in \rmin$, where $n$ is equal to the degree of $P$,
$\deg P$.
The numbers $c_1,\ldots,c_n$, 
are called the \new{corners} of $P$. 
They coincide with the points of nondifferentiability of $\hat P$.
If $c$ is a corner, the \new{multiplicity} of $c$,
which is equal to the number of indices $i$ for which
$c_i=c$, coincides, when $c\neq \zero$, with the variation of slope
of $\hat P$ at $c$, $\hat P'(c^-)-\hat P'(c^+)$.
The multiplicity of the corner $\zero$
is equal to the valuation of $P$, $\val P$.
We denote by $\gamma_1,\ldots,\gamma_N$ the corners
of the characteristic polynomial $P_A$, that we call
the (algebraic) {\em eigenvalues} of $A$. 
Note that the valuation $\val P_A$ can be computed
by introducing the matrix $\val A\in \rmin^{n\times n}$,
such that $(\val A)_{ij}=\val A_{ij}$. Then, $\val P_A$
is equal to the min-plus permanent of the matrix $\val A$.
By symmetry, the degree $\deg P_A$ is equal to the 
max-plus permanent of the matrix
$\deg A \in \rmax^{n\times n}$,
such that $(\deg A)_{ij}=\deg A_{ij}$.
The corners $\gamma_i$ (and so, the polynomial
function $\hat P_A$) can be computed
in $O(n^4d)$ time by adapting
the method of Burkard and Butkovi\v{c}~\cite{bb02}.
(It is not known whether the sequence of coefficients
of the {\em formal} polynomial $P_A$
can be computed in polynomial time.)

For any matrix $B\in \rmin^{n\times n}$ such that $\perm B\neq\zero$,
we define the graph $\opt(B)$ as the set of arcs belonging
to optimal assignments: the nodes of $\opt(B)$ 
are $1,\ldots,n$ and there is an arc from $i$ to
$j$ if there is a permutation $\sigma$ such that
$j=\sigma(i)$ and $|\sigma|_B=\perm B$.

We shall say that two vectors
$U,V$ of dimension $n$ with entries in $\rmin\setminus\{\zero\}$
form a \new{Hungarian pair} with respect to $B$ if, for all $i,j$,
we have $B_{ij}\geq U_iV_j$, 
and $U_1\cdots U_n V_1\cdots V_n=\perm B$,
the products being understood
in the min-plus sense. Thus, $(U,V)$ coincides with the
optimal dual variable in the linear programming
formulation of the optimal assignment problem (see~\cite[2.4]{br}).
In particular, a Hungarian pair always exists if
the optimal assignment problem is feasible,
i.e., if $\perm B \neq \zero$, and it can
be computed in $O(n^3)$ time by the Hungarian
algorithm.
For any Hungarian pair $(U,V)$, we now define 
the {\em saturation graph}, 
$\sat(B,U,V)$, which has nodes $1,\ldots,n$
and an arc from $i$ to $j$ if $B_{ij}=U_iV_{j}$.

For every finite corner $\gamma$ of $P_A$, we define
the digraphs $\opt_{0}(\gamma),\ldots,\opt_{d}(\gamma)$:
$\opt_{k}(\gamma)$ has nodes $1,\ldots,n$, and
an arc from $i$ to $j$ if $(i,j)\in \opt(\hat A(\gamma))$
and $(A_k)_{ij}\gamma^k=\hat A_{ij}(\gamma)$. 
For every Hungarian pair $(U,V)$ with respect
to the matrix $\hat A(\gamma)$, we also
define the digraphs $\sat_{0}(\gamma,U,V),\ldots,\sat_{d}(\gamma,U,V)$
by replacing $\opt(\hat A(\gamma))$ by $\sat(\hat A(\gamma),U,V)$,
in the definition of $\opt_{0}(\gamma),\ldots,\opt_{d}(\gamma)$.
Finally, if $G$ is any digraph with nodes $1,\ldots,n$, 
and if $b\in \C^{n\times n}$,
we define the matrix $b^{G}$, which is such that
$(b^{G})_{ij}=b_{ij}$ if $(i,j)\in G$,
and $(b^{G})_{ij}=0$ otherwise.
In the following theorem, and in the sequel,
eigenvalues are counted with multiplicities.
\begin{theorem}
\label{th-1}
Let $\gamma$ denote any finite corner 
of the min-plus characteristic polynomial
$P_A$. For every $0\leq k\leq d$,
let $G_{k}$ be equal either to 
$\opt_{k}(\gamma)$ or $\sat_{k}(\gamma,U,V)$, for any choice
of the Hungarian pair $U,V$ with respect to $\hat A(\gamma)$.
Consider the auxiliary pencil
\[
a^{(\gamma)}:= a_0^{G_{0}}  + Xa_1^{G_{1}} +\cdots + 
X^da_d^{G_{d}} \enspace .
\]
Then, if the pencil $a^{(\gamma)}$ has $m_\gamma$
non-zero eigenvalues, $\lambda_1,\ldots,\lambda_{m_\gamma}$,
the pencil $\sA_\epsilon$ has $m_\gamma$ eigenvalues
$\sL_{\epsilon,1},\ldots,\sL_{\epsilon,m_\gamma}$
with respective equivalents of the form
\(
\sL_{\epsilon,i}\sim \lambda_i \epsilon^{\gamma}
\);
if $0$ is an eigenvalue of multiplicity $m'_\gamma$
of the pencil $a^{(\gamma)}$, the pencil
$\sA_\epsilon$ has precisely $m'_\gamma$
eigenvalues $\sL_\epsilon$ such that
$\epsilon^{-\gamma}\sL_\eps$ converges to zero,
and all the other eigenvalues $\sL_\epsilon$ of $\sA_\epsilon$
are such that the modulus of $\epsilon^{-\gamma}\sL_\eps$
converges to infinity.
Moreover, for generic values of the parameters
$(a_k)_{ij}$,
$m_\gamma$ coincides with the multiplicity
of the corner $\gamma$, and $m'_\gamma$ coincides with the sum
of multiplicities of all the corners
greater than $\gamma$.
\end{theorem}
Since the sum of the multiplicities of the corners
of $P_A$ is equal generically to the degree of $\det(\sA_\epsilon)$,
and since $\sA_\epsilon$ has a number of identically zero
eigenvalues generically 
equal to the multiplicity of $\zero$ as a corner of $P_A$,
Theorem~\ref{th-1} provides generically
asymptotic equivalents of all the eigenvalues of $\sA_\epsilon$.
\begin{example}
Consider $\sA_\epsilon=\sA_{\epsilon,0}-X\id$, and
\[
\sA_{\eps,0} = \left[\begin{array}{ccc}
b_{11} & b_{12} & b_{13}\\
b_{21} & b_{22}\eps & b_{23}\eps\\
b_{31} & b_{32}\eps & b_{33}\eps
\end{array}\right]\enspace,\;\mrm{where $b_{ij}\in \C$.}
\]
The associated min-plus matrix pencil and characteristic polynomial
function are 
\[ A= \left[\begin{array}{ccc}0\oplus X&0&0\\0&1\oplus X&1\\0&1&1\oplus X
\end{array}\right]\enspace,
\qquad 
\hat{P_A}(x)=(x\oplus 0)^2(x\oplus 1)\enspace,
\]
so that the corners are $\gamma_1=\gamma_2=0$,
with multiplicity $2$, and $\gamma_3=1$, 
with multiplicity $1$. We first consider the corner
$\gamma=0$. Then $U=V=(0,0,0)$ yields
a Hungarian pair with respect to the matrix
\[
\hat A(0)= \left[\begin{array}{ccc}0_{01}&0_0&0_0
\\ 0_0&0_1&1\\ 0_0&1& 0_1
\end{array}\right]\enspace,
\]
where we adopt the following convention
to visualize the digraphs $\sat_k(U,V)$:
an arc $(i,j)$ belongs to $\sat_k(U,V)$ if $k$ is put
as a subscript of the entry $\hat A_{ij}(0)$.
For instance, $\hat{A}_{11}(0)=0$, and $(1,1)$
belongs both to $\sat_0(U,V)$ and $\sat_1(U,V)$.
Entries without subscripts, like
$\hat A_{23}(0)=1$, correspond
to arcs which do not belong to $\sat(U,V)$.
The eigenvalues of the auxiliary pencil $a^{(0)}$ are the roots of
\[
\det \left[\begin{array}{ccc}
b_{11} -\lambda & b_{12} & b_{13}\\
b_{21} & -\lambda &0 \\
b_{31} & 0 & -\lambda
\end{array}\right] 
= \lambda(-\lambda^2+\lambda b_{11}+ b_{12}b_{21}+b_{31}b_{31}) =0
\enspace.
\]
Theorem~\ref{th-1} predicts that this equation has,
for generic values of the parameters $b_{ij}$,
two non-zero roots, $\lambda_1,\lambda_2$,
which yields two eigenvalues of $\sA_\epsilon$,
$\sL_{\epsilon,m}\sim \lambda_m\epsilon^0=\lambda_m$,
for $m=1,2$.
Consider finally the corner $\gamma=1$.
We can take $U=(0,1,1)$, $V=(-1,0,0)$,
and the previous computations become
\[
\hat A(1)= \left[\begin{array}{ccc}0&0_0&0_0
\\ 0_0&1_{01}&1_0\\
0_0&
1_0&1_{01}
\end{array}\right]\enspace,\qquad
\det
\left[\begin{array}{ccc}
0 & b_{12} & b_{13}\\
b_{21} & b_{22}-\lambda &b_{23} \\
b_{31} & b_{31} & b_{33}-\lambda 
\end{array}\right] =0 \enspace .
\]
The latest equation yields
\(
\lambda(b_{12}b_{21}+b_{13}b_{31})
+b_{12}b_{23}b_{31}+b_{13}b_{32}b_{21}
-b_{21}b_{12}b_{33}
-b_{31}b_{13}b_{22}=0\).
Theorem~\ref{th-1} predicts that this equation
has generically a unique nonzero root, $\lambda_1$, 
and that there is a branch $\sL_{\epsilon,1}\sim \lambda_1 \epsilon$.
\end{example}
As a typical application of Theorem~\ref{th-1},
let us consider the following singular perturbation
of an affine pencil,
\(
\sA_\epsilon=\epsilon X^2m + Xc +k\),
already considered in~\cite{najman}. 
For non-zero values of the entries of the matrices $m,c$,
and $k$, the associated min-plus characteristic polynomial
function is $\hat P_A(x)=(0\oplus x)^n(0\oplus 1x)^n$.
Moreover, the pencils $Xc+k$ and $Xm+c$ generically both have
$n$ finite non-zero eigenvalues, denoted by
$\lambda_1,\ldots,\lambda_n$, and $\mu_1,\ldots,\mu_n$,
respectively. 
Then, it is easy to derive from Theorem~\ref{th-1} that the pencil
$\sA_\epsilon$ has $n$ eigenvalues
$\sL_{\epsilon,i}\sim \lambda_i\epsilon^0$,
and $n$ eigenvalues $\sL_{\epsilon,i}\sim \mu_i\epsilon^{-1}$.
Consider now the following non-generic situation.
Assume that the pencil $cX+k$ is given, that it is regular,
and that its Weierstrass normal form comprises
$q_0$ Jordan blocks for the eigenvalue $0$, with respective sizes
$s_0^1,\ldots,s_0^{q_0}$, and $q_\infty$
Jordan blocks for the eigenvalue $\infty$, with respective sizes
$s_\infty^1,\ldots,s_\infty^{q_\infty}$. We set
$d_0=s_0^1+\cdots+s_0^{q_0}$, 
$d_\infty=s_\infty^1+\cdots+s_\infty^{q_\infty}$.
We also denote by $q'_0$ the number of
one dimensional Jordan blocks for
the eigenvalue $0$ of the pencil $cX+k$.
We denote by $\lambda_1,\ldots,\lambda_r$
the finite non-zero eigenvalues of $cX+k$
(of course, $r+d_0+d_\infty=n$).
We also denote by $\mu_1,\ldots,\mu_t$
the finite non-zero eigenvalues 
of the pencil $Xm+c$.
We say that an eigenvalue $\sL_\epsilon$
is \new{of order} $\epsilon^\Lambda$
if $\sL_\epsilon\sim\lambda\epsilon^\Lambda$,
for some $\lambda\in\C\setminus\{0\}$.
The following result should be compared with~\cite{najman},
where partial results are obtained in a similar situation.
\begin{corollary}\label{cor-1}
The pencil $\sA_\epsilon=\epsilon X^2m + Xc +k$ has precisely
\begin{enumerate}
\renewcommand{\theenumi}{\roman{enumi}}
\item\label{i1} $r$ eigenvalues of order $\epsilon^0$, which converge respectively to $\lambda_i$, for $i=1,\ldots,r$;
\item\label{i2} $t$ eigenvalues 
or order $\epsilon^{-1}$, which are respectively
equivalent to $\mu_i\epsilon^{-1}$, for $i=1,\ldots,t$.
\end{enumerate}
It has at least
\begin{enumerate}
\renewcommand{\theenumi}{\roman{enumi}}
\setcounter{enumi}{2}
\item\label{i5} $2q_0-q'_0$ eigenvalues identically equal to zero.
\end{enumerate}
Finally, for generic values of the parameters $m_{ij}$,
we have $t=n-q_\infty$, and the pencil $\sA_\epsilon$ has
precisely:
\begin{enumerate}
\renewcommand{\theenumi}{\roman{enumi}}
\setcounter{enumi}{3}
\item\label{i3} $s^i_\infty+1$ eigenvalues of order $\epsilon^{-1/(s^i_\infty+1)}$,
for $i=1,\ldots,q_\infty$; 
\item\label{i4} $s^i_0-2$ eigenvalues of order $\epsilon^{1/(s^i_0-2)}$,
for every $i$ such that $1\leq i\leq r$ and $s^i_0> 2$.
\end{enumerate}
\end{corollary}
Corollary~\ref{cor-1} provides, for generic values of $m$,
the leading 
exponents of all the eigenvalues of the pencil $\sA_\epsilon$.
In cases~\ref{i3}--\ref{i4}, the generic values of the leading
coefficients of the eigenvalues can be determined by formul\ae\
essentially similar to the case of~\cite{lidskii,abg04a}. This will
be detailed elsewhere.

{\em Acknowledgement}.\/ The third author thanks Jean-Jacques Loiseau
for having suggested to look for a generalization of the result of~\cite{abg04a} to matrix pencils.


\begin{thebibliography}{BCOQ92}
\providecommand{\url}[1]{\texttt{#1}}
\providecommand{\urlprefix}{URL }
  \providecommand{\doi}[1]{Eprint \href{http://dx.doi.org/#1}{doi:#1}}
\providecommand{\arxiv}[2][]{Also \href{http://www.arXiv.org/abs/#2}{arXiv:#2}}

\bibitem[ABG04]{abg04a}
M.~Akian, R.~Bapat, and S.~Gaubert.
\newblock Generic asymptotics of eigenvalues and min-plus algebra.
\newblock Rapport de recherche 5104, INRIA, Le Chesnay, France, Feb. 2004.
\newblock \arxiv{math.SP/0402090}.

\bibitem[BB03]{bb02}
R.~E. Burkard and P.~Butkovi{\v{c}}.
\newblock Finding all essential terms of a characteristic maxpolynomial.
\newblock \emph{Discrete Appl. Math.}, 130(3):367--380, 2003.
\newblock \doi{10.1016/S0166-218X(03)00223-3}.

\bibitem[BCOQ92]{bcoq}
F.~Baccelli, G.~Cohen, G.~Olsder, and J.~Quadrat.
\newblock \emph{Synchronization and Linearity --- an Algebra for Discrete Event
  Systems}.
\newblock Wiley, 1992.

\bibitem[BR97]{br}
R.~Bapat and T.~Raghavan.
\newblock \emph{Nonnegative Matrices and Application}.
\newblock Cambridge University Press, 1997.

\bibitem[CGM80]{cuning80}
R.~Cuninghame-Green and P.~Meijer.
\newblock An algebra for piecewise-linear minimax problems.
\newblock \emph{Dicrete Appl. Math}, 2:267--294, 1980.

\bibitem[Lid65]{lidskii}
V.~Lidski\u\i.
\newblock Perturbation theory of non-conjugate operators.
\newblock \emph{U.S.S.R. Comput. Math. and Math. Phys.,}, 1:73--85, 1965.
\newblock (\u Z. Vy\v cisl. Mat. i Mat. Fiz. 6, no. 1, 52--60, 1965).

\bibitem[MBO97]{mbo97}
J.~Moro, J.~V. Burke, and M.~L. Overton.
\newblock On the {L}idskii-{V}ishik-{L}yusternik perturbation theory for
  eigenvalues of matrices with arbitrary {J}ordan structure.
\newblock \emph{SIAM J. Matrix Anal. Appl.}, 18(4):793--817, 1997.
\newblock \doi{10.1137/S0895479895294666}.

\bibitem[ME98]{maedelman}
Y.~Ma and A.~Edelman.
\newblock Nongeneric eigenvalue perturbations of {J}ordan blocks.
\newblock \emph{Linear Algebra Appl.}, 273:45--63, 1998.
\newblock \doi{10.1016/S0024-3795(97)00342-X}.

\bibitem[Naj99]{najman}
B.~Najman.
\newblock The asymptotic behavior of the eigenvalues of a singularly perturbed
  linear pencil.
\newblock \emph{SIAM J. Matrix Anal. Appl.}, 20(2):420--427, 1999.
\newblock \doi{10.1137/S0895479896299949}.

\bibitem[VL60]{vishik}
M.~I. Vi{\v{s}}ik and L.~A. Ljusternik.
\newblock Solution of some perturbation problems in the case of matrices and
  self-adjoint or non-selfadjoint differential equations. {I}.
\newblock \emph{Russian Math. Surveys}, 15(3):1--73, 1960.

\end{thebibliography}
\end{document}